\documentclass[11pt]{article}

\usepackage{amsmath,amssymb,amsthm,amsfonts,verbatim}
\usepackage{graphics,pst-plot,psfrag,epsfig}
\usepackage{graphicx}
\usepackage[all,2cell,dvips]{xy}

\title{The lower central series and pseudo-Anosov dilatations}
\author{Benson Farb, Christopher J. Leininger, and Dan Margalit
\thanks{The authors gratefully acknowledge support from the National
Science Foundation.}}

\theoremstyle{plain}
\newtheorem{theorem}{Theorem}[section]

\newtheorem{proposition}[theorem]{Proposition}

\newtheorem{lemma}[theorem]{Lemma}

\newtheorem{conjecture}[theorem]{Conjecture}
\newtheorem{question}[theorem]{Question}

\newtheorem*{theorem:pc}{Theorem \ref{theorem:principal congruence}}
\newtheorem*{theorem:jk}{Theorem \ref{theorem:johnson}}
\newtheorem*{theorem:br}{Theorem \ref{theorem:brunnian}}

\newcommand{\nc}{\newcommand}
\nc{\dmo}{\DeclareMathOperator}

\nc{\C}{\mathcal{C}}
\nc{\I}{\mathcal{I}}
\nc{\K}{\mathcal{K}}
\nc{\N}{\mathcal{N}}
\nc{\R}{\mathbb{R}}
\nc{\Z}{\mathbb{Z}}
\dmo{\PSL}{PSL}
\dmo{\Teich}{Teich}
\dmo{\spec}{spec}
\nc{\gin}{i}
\nc{\ga}{\Gamma}
\dmo{\Out}{Out}
\dmo{\brun}{Brun}

\dmo\Sp{Sp}
\dmo\Mod{Mod}
\dmo\PMod{PMod}
\dmo\genus{genus}
\dmo\htop{h_{top}}
\dmo\CAT{CAT}

\def\remark{{\bf {\bigskip}{\noindent}Remark. }}

\nc{\margin}[1]{\marginpar{\scriptsize #1}}

\begin{document}
\maketitle
\begin{abstract}
The theme of this paper is that algebraic complexity implies dynamical complexity for pseudo-Anosov homeomorphisms of a
closed surface $S_g$ of genus $g$.  Penner proved that the logarithm of the minimal dilatation for a pseudo-Anosov
homeomorphism of $S_g$ tends to zero at the rate $1/g$. We consider here the smallest dilatation of any pseudo-Anosov
homeomorphism of $S_g$ acting trivially on $\Gamma/\Gamma_k$, the quotient of $\Gamma = \pi_1(S_g)$ by the
$k^{\mbox{\tiny th}}$ term of its lower central series, $k \geq 1$.  In contrast to Penner's asymptotics, we prove that
this minimal dilatation is bounded above and below, independently of $g$, with bounds tending to infinity with $k$. For
example, in the case of the Torelli group $\I(S_g)$, we prove that $L(\I(S_g))$, the logarithm of the minimal
dilatation in $\I(S_g)$, satisfies $.197 < L(\I(S_g))< 4.127$. In contrast, we find pseudo-Anosov mapping classes
acting trivially on $\Gamma/\Gamma_k$ whose asymptotic translation lengths on the complex of curves tend to $0$ as
$g\to\infty$.
\end{abstract}

\maketitle


\section{Introduction}

Let $\Mod(S)$ denote the mapping class group of a closed, orientable surface $S = S_g$ of genus $g\geq 2$; this is the
group of homotopy classes of orientation preserving homeomorphisms of $S$. According to the Nielsen--Thurston
classification, every mapping class $f\in \Mod(S)$ which is not finite order and is not reducible (i.e. does not fix
the isotopy class of any essential 1-submanifold) is {\em pseudo-Anosov}, i.e. it has a representative which is a pseudo-Anosov
homeomorphism; see \cite{FLP,Th}.

Attached to each pseudo-Anosov $f \in \Mod(S)$ is its {\em dilatation} $\lambda(f)$. This is an algebraic integer which
records the exponential growth rate of lengths of curves under iteration of $f$, in any fixed metric on $S$; see
\cite{Th}.  The number $\log(\lambda(f))$ equals the minimal topological entropy of any element in the homotopy class
$f$; this minimum is realized by a pseudo-Anosov homeomorphism representing $f$ (see \cite[Expos\'e 10]{FLP}). From
another perspective, $\log(\lambda(f))$ is the translation length of $f$ as an isometry of the {\em Teichm\"uller space
of $S$} equipped with the Teichm\"uller metric.

Following Penner, we consider the set
$$\spec(\Mod(S))=\{\log(\lambda(f)): \mbox{$f\in \Mod(S)$ is
pseudo-Anosov}\} \subset \R$$ which can be thought of as the {\em length
spectrum} of the moduli space of genus $g$ Riemann surfaces.  We will
also consider, for various subgroups $H<\Mod(S)$, the subset $\spec(H)
\subset
\spec(\Mod(S))$ obtained by restricting to pseudo-Anosov elements of
$H$.  Arnoux--Yoccoz \cite{AY} and Ivanov
\cite{Iv1} proved that $\spec(\Mod(S))$ is discrete as a subset of $\R$.  It follows that for any subgroup $H<\Mod(S)$,
the set $\spec(H)$ is empty or has a least element $L(H)$.

If $F(x)$ and $G(x)$ are any real-valued functions, we write $F(x)
\asymp G(x)$ if there exists a $C>0$ so that
$F(x)/G(x)\in [1/C,C]$ for all $x$.  Penner \cite{Pe} proved that
$$L(\Mod(S_g)) \asymp \frac{1}{g}.$$
\noindent In particular, as genus increases, there are pseudo-Anosov
mapping classes with dilatations arbitrarily close to 1.

\paragraph{Torelli dilatations. }
The theme of this paper is that algebraic complexity implies dynamical complexity for pseudo-Anosov homeomorphisms. The
following contrast to Penner's theorem is a first instance of this phenomenon.  Below, $\I(S)$ denotes the
\emph{Torelli group}, which is defined to be the subgroup of $\Mod(S)$ consisting of elements which act trivially on
$H_1(S;\Z)$.

\begin{theorem}
\label{theorem:torelli} For $g \geq 2$, we have
$$ .197 < L(\I(S_g)) < 4.127.$$
\end{theorem}

The main point of Theorem \ref{theorem:torelli} is that $L(\I(S_g)) \asymp 1$; in other words, the bounds given in
Theorem~\ref{theorem:torelli} are universal with respect to $g$.  In contrast, Theorem \ref{theorem:LCC curve complex}
below states that the minimal translation length in the complex of curves for pseudo-Anosov mapping classes in
$\I(S_g)$ tends to $0$ as $g \to\infty$.

We remark that every pseudo-Anosov element $f \in \I(S)$ has nonorientable stable and unstable foliations since
otherwise $\lambda(f)$ would be a nontrivial eigenvalue for the action on homology; see \cite{Th}.  However, this
condition alone is insufficient to guarantee uniform upper and lower bounds for $\log(\lambda(f))$.  For example, a
construction of McMullen \cite{Mc} can be used to produce a sequence of
pseudo-Anosov elements $f_n \in \Mod(S_{g_n})$, where $g_n \to
\infty$, each
$f_n$ has nonorientable foliations, and $\log(\lambda(f_n)) \asymp
1/g_n$.

For the \emph{Johnson kernel}, which is the subgroup $\K(S)$ of $\I(S)$ generated by Dehn twists about separating curves, we obtain slightly better bounds, $.693 < L(\K(S)) < 4.127$; see Proposition~\ref{proposition:johnson} below.

\paragraph{The Johnson filtration. }
The groups $\I(S)$ and
$\K(S)$ are the first terms of the \emph{Johnson filtration} of
$\Mod(S)$, which is the sequence of groups
\[ \N_k(S) = {\rm kernel}(\Mod(S)\to \Out(\Gamma/\Gamma_k)) \]
where $\Gamma_k$ is the $k^{\mbox{\tiny th}}$ term of the lower central
series for $\Gamma = \pi_1(S)$, defined inductively by $\Gamma_0=\Gamma$
and $\Gamma_{k+1}=[\Gamma_k,\Gamma]$.  It is a classical theorem of
Magnus that $\{\Gamma_k\}$ is a \emph{filtration} of $\Gamma$, which
means that $\Gamma_{k+1} < \Gamma_k$ and $\bigcap_{k=1}^\infty
\Gamma_k = 1$; it follows that $\{\N_k(S)\}$ is a filtration of $\Mod(S)$.  By definition, $\N_0(S)=\Mod(S)$
and $\N_1(S)=\I(S)$.  It is a theorem of Johnson \cite{Jo2} that $\N_2(S)$ is isomorphic to $\K(S)$. It is a fact that
$\{\N_k(S)\}$ is a \emph{central} filtration of $\Mod(S)$ (i.e. successive quotients are abelian) \cite{BL}, and so
$\N_{k+1}(S)$ contains the $k^{\mbox{\tiny th}}$ term of the lower central series of the Torelli group $\I(S)$ (the
lower central series descends faster than any central series).

For a fixed surface $S$, a compactness argument (see Proposition \ref{proposition:filtration} below) readily gives that
$L(\N_k(S))\to \infty$ as $k\to \infty$; that is, as one specifies more and more algebraic conditions by considering
pseudo-Anosov homeomorphisms fixing deeper quotients $\Gamma/\Gamma_k$, the corresponding dynamical complexity (measured as the
dilatation) must diverge to infinity. Our main result is that this divergence is uniform over all surfaces.

\begin{theorem}
\label{theorem:johnson filtration}
Given $k \geq 1$, there exist $M(k)$ and $m(k)$, where $m(k) \to
\infty$ as $k \to \infty$, so that
\[ m(k) < L(\N_k(S_g))  < M(k) \]
for every $g \geq 2$.
\end{theorem}

Again, we compare Theorem \ref{theorem:johnson filtration} with Theorem~\ref{theorem:LCC curve complex} below.

We do not have good control over the constants $m(k)$ and $M(k)$ in this theorem, and we are interested in more precise
asymptotics for $L(\N_k(S_g))$ as $g \to \infty$.  We pose the following.

\begin{question}
Let $k$ be fixed.  As we increase the genus $g$, what is $\inf\,L(\N_k(S_g))$? What is $\sup\,L(\N_k(S_g))$?  Does
$\lim\,L(\N_k(S_g))$ exist? Are any of these quantities realized for some $g$?  Of particular interest is $L(\I(S_g))$.
\end{question}

\medskip

We consider Theorem~\ref{theorem:torelli} (and Proposition~\ref{proposition:johnson} below) as warmups for Theorem
\ref{theorem:johnson filtration}, as their proofs contain many of the main ideas.  Moreover, in these cases we compute
explicit bounds, whereas for arbitrary values of $k$ we do not.

The upper bound for $L(\I(S))$ and $L(\K(S))$ in Theorem~\ref{theorem:torelli} and Proposition
\ref{proposition:johnson} is given by explicit construction; see \S\ref{section:upper} below.   In addition, this
construction is used to derive the upper bound in Theorem \ref{theorem:johnson filtration}, using the relationship
between the lower central series of $\I(S)$ and $\{ \N_k(S)\}$. One easily checks that this upper bound grows at most
exponentially with $k$; see \S\ref{section:asymptotic upper bounds}. The proof of the lower bound begins with the
following.

\begin{proposition} \label{proposition:torelli intersections}
Suppose $f\in \I(S)$ is pseudo-Anosov.  If $c$ is a separating curve, then $\gin(c,f(c)) \geq 4$.  If $c$ is a
nonseparating curve, then $\gin(c,f^j(c)) \geq 2$ for $j=1$ or $j=2$.
\end{proposition}

The idea is to use this proposition, combined with a surgery argument, to find a curve whose length in a certain metric
is stretched by a definite amount under a pseudo-Anosov mapping class.  The metric comes from a quadratic differential
with vertical and horizontal foliations given by the stable and unstable foliations for the mapping class. The
relationship between the metric and the foliations implies that the amount of stretching bounds the dilatation from
below; see Lemma \ref{lemma:wolpert lemma}.

The proof of the lower bound in Theorem \ref{theorem:johnson filtration} follows a similar line of reasoning and
requires an asymptotic version of Proposition \ref{proposition:torelli intersections}. We show that, for $f$ lying in a
deep term of the Johnson filtration, $\gin(c,f(c))$ is large for every curve $c$ with $f(c) \neq c$
(Lemma~\ref{lemma:asymptotic intersections}).

\paragraph{Translation lengths on the complex of curves.}

One can also consider the global topological complexity of a pseudo-Anosov homeomorphism given by the translation
lengths on the ($1$-skeleton of the) {\em complex of curves} $\C=\C(S)$. This complex, defined by Harvey \cite{H}, has
a vertex for each isotopy class of essential simple closed curves in $S$ and an edge for each pair of vertices with
disjoint representatives.  We endow $\C$ with the path metric $d_{\C}$ (after declaring each edge to have length 1) and
define the \emph{asymptotic translation length} for the action of $f$ on $\C$ by
\[
\tau_\C(f) = \liminf_{j \to \infty} \frac{d_{\C}(c,f^j(c))}{j}
\]
for any curve $c$ (this is independent of the choice of curve $c$).  For any subgroup $H<\Mod(S)$, we denote by
$L_\C(H)$ the infimum of $\tau_\C(f)$ over all pseudo-Anosov elements $f \in H$ . Masur--Minsky \cite[Prop 4.6]{MM} proved that for any fixed $g$, $L_\C(\Mod(S_g))>0$.

Our first result in this direction shows that
$L_\C(\Mod(S_g))$ tends to $0$ strictly faster than
$L(\Mod(S_g))\asymp 1/g$.

\begin{theorem}
\label{theorem:curve complex infs}
For any $g \geq 2$, we have
$$L_\C(\Mod(S_g)) < \frac{4 \log(2+\sqrt 3)}{g \log\left(g-\frac{1}{2}\right)}.$$
\end{theorem}

The following result provides a contrast to
Theorem \ref{theorem:torelli} and Theorem
\ref{theorem:johnson filtration}.

\begin{theorem}
\label{theorem:LCC curve complex}
For any $k$, we have
$$L_\C(\N_k(S_g)) \to 0$$
as $g \to \infty$.
\end{theorem}

\paragraph{Congruence subgroups. }
The ideas involved in the proof of Theorem~\ref{theorem:johnson filtration} provide bounds for a different sequence of
subgroups of $\Mod(S)$.  Let $\Mod(S)[r]$ denote the {\em principal level $r$ congruence subgroup} of $\Mod(S)$, which
is defined to be the finite index subgroup of $\Mod(S)$ consisting of those elements acting trivially on
$H_1(S;\Z/r\Z)$.  We prove the following in \S\ref{section:principal congruence}.

\begin{theorem}
\label{theorem:principal congruence} If $g \geq 2$ and $r \geq 3$, then $$.197< L(\Mod(S_g)[r])< 4.127.$$
\end{theorem}

Theorem \ref{theorem:principal congruence} puts strong constraints
on the possibilities for pseudo-Anosov elements of least dilatation in
$\Mod(S)$.

\medskip
\noindent {\bf Brunnian subgroups.} In \S\ref{section:brunnian} we provide a different illustration of our theme by
considering pseudo-Anosov mapping classes in the {\em Brunnian subgroup} $\brun(S_{g,p})$ of the mapping class group of the orientable
surface $S_{g,p}$ of genus $g$ with $p>0$ punctures. This is the subgroup consisting of those mapping classes which are
isotopic to the identity once any puncture is filled in (see \S\ref{section:brunnian} for details).

\begin{theorem}
\label{theorem:brunnian} For any $g\geq 0$ and any $p\geq 5$, we have
$$L(\brun(S_{g,p})) > \log \left(\frac{p}{4}\right).$$
\end{theorem}

\paragraph{Related results in the literature. }
As we noted above, Penner \cite{Pe} gave the first proof that $L(\Mod(S_g)) \asymp 1/g$.  His upper bound was improved
upon by Bauer \cite{Ba1,Ba2}, who gave new examples with small dilatation.  McMullen \cite{Mc} gave a different
construction for the upper bound of Penner's asymptotics using fibered $3$-manifolds with infinitely many fibrations.
Brinkmann \cite{Br}, Hironaka--Kin \cite{HK}, and Minakawa \cite{Mk} also gave examples proving the same upper bound
for the asymptotics. The best known general upper bound is $\log(2 + \sqrt 3)/g$ given by Hironaka--Kin \cite{HK} and
Minakawa \cite{Mk}.  The precise value of $L(\Mod(S_g))$ is not known for any $g>1$; some related values have been
calculated \cite{Zh,SKL,HS}.

The second author \cite{Le} investigated the question of the minimal dilatation for the class of subgroups $\langle
T_A,T_B \rangle$ generated by two positive multitwists $T_A$ and $T_B$.  In this case the infimum of $L(\langle T_A,T_B
\rangle)$ over all genus and all such subgroups is the logarithm of Lehmer's number $\log(\lambda_L) \approx .162$, and
is realized on a genus $5$ surface.  For pure braid groups $PB_n$, Song \cite{So} proved that $\log(2+\sqrt{5}) \leq
L(PB_n)$.  In fact, one has $L(PB_n) \asymp 1$; see \S\ref{section:upper} for an upper bound. Finally, for the
hyperelliptic subgroups, Hironaka--Kin \cite{HK} proved that the asymptotics are the same as those of $\Mod(S_g)$,
giving an explicit upper bound of $\log(2 + \sqrt{3})/g$. Moreover, their examples descend to braids that cyclically
permute the punctures. Thus, they also obtain an upper bound $L(PB_{2g+1}) \leq (2+1/g)\log(2+\sqrt{3}) \leq 5
\log(2+\sqrt{3})/2$.

\paragraph{Acknowledgements.} The authors are grateful to Justin
Malestein for reading an earlier draft and making corrections.  We also
thank John Franks for a suggestion which improved our lower bounds by a
factor of two; see the remark at the end of \S \ref{section:proof of
lower bound}.  Finally, we thank the referee for numerous comments and
suggestions which greatly improved the final version of this paper.

\section{Torelli groups}
\label{section:torelli groups}

The mechanism by which we prove that a pseudo-Anosov element $f\in\I(S)$ is forced to have big
dilatation comes from the action of $f$ on simple closed curves.  This
is best explored via intersection numbers.

\subsection{\boldmath$\I(S)$ and geometric intersection numbers}
\label{section:geometric intersection numbers}

Let $a$ and $b$ be free homotopy classes of simple closed curves in $S$.  The {\em geometric intersection number}
$\gin(a,b)$ is defined by
$$\gin(a,b)=\min\{|\alpha \cap \beta|: \alpha \in a \ \mbox{and\ }
\beta \in b\} .$$ We generally do not distinguish between homotopy classes of simple closed curves and particular
representatives of the classes, referring to both simply as ``curves'', with usage dictating what is meant (likewise
for mapping classes and representative homeomorphisms).  Representative
curves $a$ and $b$ of homotopy classes of the same names are in
\emph{minimal position} if they are transverse and $\gin(a,b) = |a \cap
b|$.  Whenever considering representatives of a pair of homotopy classes
we assume that they are in
minimal position unless stated otherwise.

The mod 2 intersection number of two curves can be computed as geometric intersection number mod 2 or algebraic intersection number mod 2.
Therefore we have the following fact which we will use repeatedly without mention.

\begin{center}
\textit{The algebraic and geometric intersection number of a pair of curves have the same parity.}
\end{center}

\noindent
In particular if $a$ and $b$ are homologous, or if $b$ is separating, then $i(a,b)$ is even.

We will also need the following general fact about geometric intersection
numbers.  We say that a collection of curves \emph{fills} a closed surface if
the complement is a disjoint union of disks.

\begin{lemma}
\label{lemma:filling curves intersect} If $a$ and $b$ are two simple closed curves which together fill the closed
surface $S_g$, then $\gin(a,b) \geq 2g-1$.  More generally, for any two curves $a$ and $b$ in $S$, $\gin(a,b) =
-\chi(N)$, where $N$ is a regular neighborhood of $a \cup b$.
\end{lemma}

\begin{proof}
If $a$ and $b$ fill $S_g$, then $\chi(N) < \chi(S_g) = 2-2g$ since
$S_g$ is obtained from $N$ by gluing disks to the boundary components.
Because these numbers are integers, the first statement follows from
the second.

To prove the second statement, we simply note that $N$ deformation
retracts onto $a \cup b$, thought of as a graph in $N$.  Then $\chi(N)
= V -E$, where $V = \gin(a,b)$ is the number of vertices and $E$ is
the number of edges.  Because each vertex of $a \cup b$ is $4$-valent
we see that $E = 2V = 2\gin(a,b)$, and hence $\chi(N) = -\gin(a,b)$.
\end{proof}

\begin{lemma}
\label{lemma:separating curves intersect image}
Suppose that $f \in \I(S)$, that $c$ is a separating curve, and that
$f(c)\neq c$.  Then $\gin(f(c),c)\geq 4$.
\end{lemma}

This lemma is sharp when the genus of $S$ is at least 4, for in this case one can find two separating curves $c$ and
$d$ with $\gin(c,d)=2$.  Since $d$ is trivial in homology, the Dehn
twist $T_d$ about $d$ is in $\I(S)$.  Drawing the picture, we find
$\gin(T_d(c),c) = 4$.

\begin{figure}[ht]
\begin{center}
\psfrag{c}{$c$}
\psfrag{fc}{$f(c)$}
\psfrag{b1}{$b_1$}
\psfrag{b2}{$b_2$}
\psfrag{b3}{$b_3$}
\psfrag{b4}{$b_4$}
\includegraphics[scale=.4]{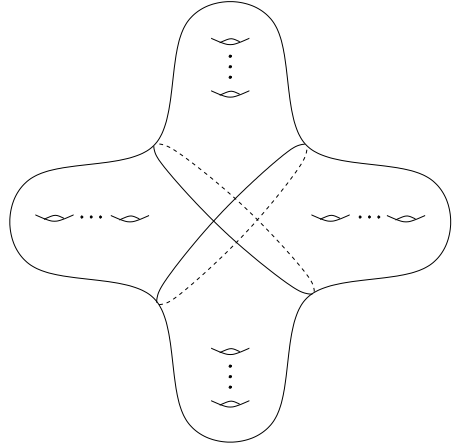}
\caption{A pair of separating curves which intersect twice.}
\label{figure:separating proof}
\end{center}
\end{figure}

\begin{proof}
First, any two separating curves have even geometric intersection number. Two distinct separating curves with intersection number zero clearly induce different splittings
of $H_1(S;\Z)$, so it is impossible to have $\gin(c,f(c))=0$.  Any two separating curves with intersection number 2
are, after composing with a homeomorphism, given as in Figure~\ref{figure:separating proof}. Again we see that they
induce different homology splittings, and so we cannot have $\gin(c,f(c))=2$.
\end{proof}

\begin{lemma}
\label{lemma:nonseparating curves intersect image}
Suppose that $f \in \I(S)$, that $c$ is a nonseparating curve, and that $f(c)\neq c$.  Then at least one of
$\gin(f(c),c)$ and $\gin(f^2(c),c)$ is at least 2.
\end{lemma}

Note that this lemma is sharp for $g \geq 3$ in the sense that there
exists an element $f \in \I(S)$ and a curve $c$ so that $\gin(c,f(c))=0$
and $\gin(c,f^2(c))=2$.  Consider, for instance, a {\em bounding pair}
$\{d,e\}$, i.e. a pair of disjoint, nonhomotopic, homologous,
nonseparating simple closed curves, and a curve $c$ that intersects both
$d$ and $e$ exactly once each; in this case, $\gin(c,T_dT_e^{-1}(c))=0$
and $\gin(c,(T_dT_e^{-1})^2(c))=2$.

\begin{figure}[ht]
\begin{center}
\psfrag{c}{$c$}
\psfrag{fc}{$f(c)$}
\psfrag{ffc}{$f^2(c)$}
\includegraphics[scale=.4]{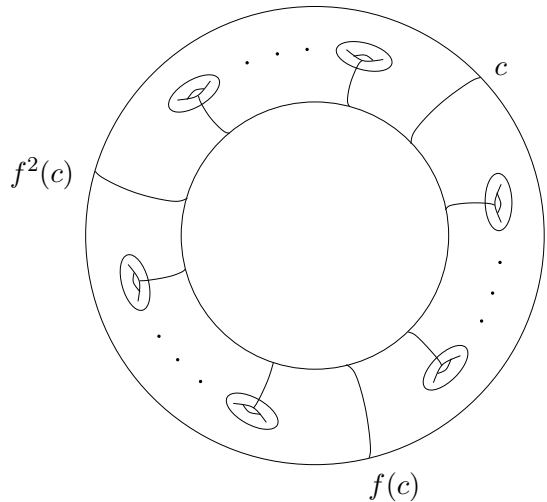}
\caption{The picture for Lemma~\ref{lemma:nonseparating curves
    intersect image}.  The unlabeled curves are the $u_i$ and $v_i$.}
\label{figure:nonseparating proof}
\end{center}
\end{figure}

\begin{proof}
Since $f\in\I(S)$ and $c \neq f(c)$ we know that $c \neq f^2(c)$
(combine Theorem 3 and Corollary 3.7 of
\cite{Iv2}). Now suppose $\gin(f(c),c)=0$ and $\gin(f^2(c),c)=0$, so that $\{c,f(c),f^2(c)\}$ is a collection of 3
distinct, disjoint simple closed curves, each representing a fixed
nonzero element of $H_1(S;\Z)$.
We can choose curves $u_2, v_2, \dots, u_g,v_g$ that are each disjoint
from $c$, $f(c)$, and $f^2(c)$, and whose corresponding homology classes
span a codimension 2 subspace $V$ of $H_1(S;\Z)$.  Now, $f$ takes the
pair $\{c,f(c)\}$ to the pair $\{f(c),f^2(c)\}$.  However, it is clear
that these two pairs induce different splittings of $V$, and so we have
a contradiction (see Figure~\ref{figure:nonseparating proof}).  Since
$f\in \I(S)$, we cannot have $\gin(f(c),c)=1$ or $\gin(f^2(c),c)=1$ since these intersection numbers must be even, so
we are done.
\end{proof}

Lemma \ref{lemma:separating curves intersect image} and Lemma
\ref{lemma:nonseparating curves intersect image} together
prove Proposition \ref{proposition:torelli
  intersections}, since a pseudo-Anosov mapping class does not fix any
  curve in the surface.

\begin{figure}[ht]
\begin{center}
\psfrag{D}{$d$} \psfrag{D'}{$d'$} \psfrag{E}{$e$} \psfrag{E'}{$e'$} \psfrag{c}{$c$} \psfrag{a}{$a$} \psfrag{a'}{$a'$}
\psfrag{b}{$b$} \psfrag{b'}{$b'$} \psfrag{fc}{$c'$}
\includegraphics[scale=.6]{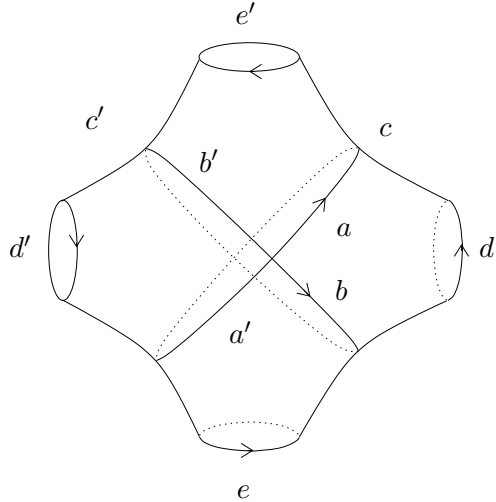}
\caption{Homologous curves $c$ and $c'$ with $\gin(c,c') = 2$ and the 4 holed sphere (the labels $a$, $a'$, $b$, and $b'$ are used in the proof of Proposition~\ref{proposition:torelli lower bound}).} \label{figure:lantern}
\end{center}
\end{figure}

We require one final fact regarding geometric intersection numbers.  Suppose $c$ and $c'$ are homologous
curves with $i(c,c') = 2$.  It follows that the two points of intersection must have opposite signs, and a regular neighborhood of $c \cup c'$ is a 4
holed sphere.  We orient $c$ and $c'$ so that $[c] = [c']$ in $H_1(S;\mathbb Z)$.  Label and orient the four boundary
components of the 4 holed sphere $d$, $d'$, $e$ and $e'$ as in Figure \ref{figure:lantern}.
\begin{lemma} \label{lemma:homologous curves git2}
Suppose $c$ and $c'$ are homologous nonseparating curves with $i(c,c') =
2$.  Suppose that $d, d' e$ and $e'$ are the boundary components of 
a $4$-holed sphere as shown in Figure~\ref{figure:lantern}.  
Then $d$ and $d'$ are
separating in $S$, and $[e] = - [e'] = [c] = [c']$ in $H_1(S;\mathbb
Z)$.
\end{lemma}
\begin{proof}
There are two pairs of pants in the 4 holed sphere which determine relations $[c] + [e'] + [d'] = 0$ and $[c'] + [e'] +
[d] = 0$.  It follows that $[d] = [d']$. Therefore $d \cup d'$ is the boundary of two subsurfaces.  The subsurface containing the
4 holed sphere is to the left of both $d$ and $d'$, so $[d] + [d'] = 0$, and hence $[d] = [d'] = 0$. A
third pair of pants defines the relation $[c'] - [e] - [d'] = 0$, and so $[e] - [e'] = [c] + [c'] = 2[c] \neq 0$. Since
the 4 holed sphere defines the relation $0 = [e] + [e'] + [d] + [d']$ and 
$[d] = [d'] = 0$, we see that $[e] = -[e'] = [c] = [c']$.
\end{proof}

\subsection{Proof of the lower bound}
\label{section:proof of lower bound}

We begin by recalling a few definitions and facts; see \cite{Ab} for a more detailed discussion.  If $f$ is
pseudo-Anosov, we let $q = q_f$ denote a holomorphic quadratic differential for which the vertical and horizontal
foliations are precisely the stable and unstable foliations for $f$, respectively. The differential $q$ determines a
euclidean cone metric, which we also denote $q$, and $f$ acts as an affine diffeomorphism (off the singularities) whose
derivative has eigenvalues $\lambda(f)$ and $\lambda(f)^{-1}$.

For any curve $c$ in $S$, we let $\ell_q(c)$ denote the infimum of $q$-lengths of representatives of $c$, which is
equivalently the length of a $q$-geodesic representative for $c$.  We note that, in general, the $q$-geodesic representative of a simple closed curve need not be embedded.


\begin{lemma}
\label{lemma:wolpert lemma} Let $f\in \Mod(S)$ be pseudo-Anosov and let $q = q_f$.  Then for any closed curve $c$ in
$S$, we have
\[ \frac{\ell_q(f(c))}{\ell_q(c)} < \lambda(f). \]
\end{lemma}

\begin{proof}
Note that because $f$ is affine with respect to $q$, the image of a geodesic representative for $c$ is a geodesic
representative for $f(c)$. Furthermore, since the leading eigenvalue of the derivative of $f$ is $\lambda =
\lambda(f)$, the length of the curve $f(c)$ differs from that of $c$ by at most a factor of $\lambda$.  Moreover, only
geodesics which are everywhere tangent to the eigenspace for $\lambda$ can be maximally stretched.  However any such
geodesic is a leaf of the stable or unstable foliation, and hence cannot be part of a closed geodesic, so the
inequality is strict.
\end{proof}

We are now ready to give the proof of the lower bound on $L(\I(S))$
given in Theorem~\ref{theorem:torelli}.

\begin{proposition}
\label{proposition:torelli lower bound} If $g \geq 2$, then $L(\I(S_g)) > .197$.
\end{proposition}

\begin{proof}
Let $f$ be an arbitrary pseudo-Anosov element of $\I(S)$.  Let $q =
q_f$, and let $c$ be a shortest curve in $S$ with respect to $q$.  We
will assume in what follows that all closed $q$-geodesics under
consideration are embedded and that all pairs of $q$-geodesics are in
minimal position. This is not true in general, but so as not to
disrupt the flow of ideas we make this assumption.  We will discuss
the minor modifications needed for the general case at the end of the
proof.

\bigskip
\noindent
{\bf Case 1. } $\gin(c,f(c))\geq 4$ or $\gin(c,f^2(c))\geq 4$.

\medskip
Let $h$ be either $f$ or $f^2$, where $\gin(c,h(c)) \geq 4$. The intersection points $c \cap h(c)$ cut each of $c$ and
$h(c)$ into arcs.  Since there are at least 4 intersection points, there is an arc $a$ of $h(c)$ which satisfies
\[
\ell_q(a) \leq \frac{\ell_q(h(c))}{4} < \frac{\lambda(h)\ell_q(c)}{4}
\]
where the second inequality comes from an application of Lemma~\ref{lemma:wolpert lemma}.  Here we have written
$\ell_q(a)$ to denote the $q$-length of the segment $a$. The endpoints of $a$ cut $c$ into two arcs.  One of which,
call it $b$, has length at most $\ell_q(c)/2$. The union $a \cup b$ is a simple closed curve in $S$.  It is nontrivial
for otherwise it would bound a disk, which we could use as a homotopy to show $\gin(c,h(c))< |c \cap h(c)|$. Since $c$
is a shortest curve with respect to $q$, we have
\begin{equation*}
\begin{array}{rcl} \ell_q(c) \leq \ell_q(a \cup b) &\leq& \ell_q(a) + \ell_q(b) \\
&\\
&<& \displaystyle \frac{\lambda(h)\ell_q(c)}{4} + \frac{\ell_q(c)}{2}\\
&\\
&=&
  \ell_q(c)\left ( \displaystyle \frac{\lambda(h)}{4} + \frac{1}{2}
\right ). \end{array}
\end{equation*}
It follows that
\[ \frac{\lambda(h)}{4} + \frac{1}{2} > 1 \]
and so $\lambda(h) > 2$.  Since $h$ is $f$ or $f^2$ and since $\lambda(f^2)=\lambda(f)^2$, we have $\lambda(f) >
\sqrt{2}$.

\bigskip
\noindent
{\bf Case 2. }$c$ is nonseparating and $\gin(c,f(c))$ and
$\gin(c,f^2(c))$ are both less than 4.

\medskip
By Lemma~\ref{lemma:nonseparating curves intersect image}, either $\gin(c,f(c))=2$ or $\gin(c,f^2(c))=2$.  Let $h$ be
either $f$ or $f^2$, where $\gin(c,h(c)) = 2$.

Let $d$ and $d'$ be the separating curves from Lemma \ref{lemma:homologous curves git2} with $c' = h(c)$.
Alternatively, the intersection points $c \cap h(c)$ define 2 arcs of $c$, say $a$ and $a'$, and two arcs of $h(c)$,
say $b$ and $b'$ as in Figure~\ref{figure:lantern}.  The curves $d$ and $d'$ are then $d = a \cup b$ and $d' =
a' \cup b'$.

Since
\[
\begin{array}{rcl} \ell_q(a) + \ell_q(b) + \ell_q(a') + \ell_q(b') &=&
 \ell_q(c) + \ell_q(h(c)) \\ & < & \ell_q(c) +
\lambda(h)\ell_q(c) \end{array}\] it follows that at least one of $d$ and $d'$, say $d$, has length bounded above by
half of $\ell_q(c) + \ell_q(h(c))$:

\begin{equation}
\label{deltabound} \ell_q(d) < \frac{\ell_q(c) + \lambda(h)\ell_q(c)}{2} \leq \frac{\ell_q(c)
+\lambda(f)^2\ell_q(c)}{2}.
\end{equation}

We now consider $d$, which is a separating curve, and its image $f(d)$, which intersect in at least four points by
Lemma \ref{lemma:separating curves intersect image}.  As in Case 1, if $a$ is the shortest arc of $f(d)$ and $b$ is the
shortest arc of $d$ cut off by $a$, then
\begin{equation}\label{alphabetabound} \ell_q(a \cup b) <  \ell_q(d)\left (
\frac{\lambda(f)}{4} + \frac{1}{2}  \right ). \end{equation} Note that we can always use $f$ (as opposed to $f^2$) since
$d$ is separating.

Also, since $c$ is shortest, we have

\begin{equation}
\label{gammabound} \ell_q(c) \leq \ell_q(a \cup b). \end{equation} Combining (\ref{deltabound}), (\ref{alphabetabound}),
and (\ref{gammabound}) we see that
\[  \ell_q(c) <  \frac{\ell_q(c) +
\lambda(f)^2\ell_q(c)}{2}\left ( \frac{\lambda(f)}{4} + \frac{1}{2}  \right ).\] In other words,
\[ \lambda(f)^3 + 2\lambda(f)^2 + \lambda(f)-6 > 0. \]

The cubic polynomial in $\lambda(f)$ on the left has one real root, and so
\[ \lambda(f) >  -\frac{2}{3}+\frac{1}{3}
\sqrt[3]{82-9\sqrt{83}}+\frac{1}{3}\sqrt[3]{82+9\sqrt{83}}\approx 1.218
\]
approximated from below.

By Proposition~\ref{proposition:torelli intersections}, these are all cases and so, after taking the logarithms, we are
done if all $q$-geodesics are embedded and all pairs are in minimal position.

\medskip
In the general case, we approximate the $q$-metric on $S$ by a
nonpositively curved Riemannian metric $q_0$ which agrees with the
$q$-metric in the complement of a small neighborhood of the singular
points.  This can be done by an explicit computation; compare, e.g.,
\cite{BH} or \cite{GT}.  Given any positive number $R > 0$, which is not
one of the $q$-lengths of a curve, we can choose this approximation so
that the set of curves with $q$-length at most $R$ is precisely the same
as the set of curves with $q_0$-length at most $R$. Moreover, given
$\epsilon > 1$, we may assume that the ratio of $q$-length and
$q_0$-length of any curve is between $\epsilon$ and $1/\epsilon$.  In
particular, since we can assume that $\lambda(f) \leq 2$, say, then we
may choose $q_0$ so that for the finite set of curves with $q_0$-length
at most $R$, we have
\[
\frac{\ell_{q_0}(f(c))}{\ell_{q_0}(c)} < \lambda(f).
\]
Since a geodesic representative of any simple closed curve in a nonpositively curved Riemannian metric on a surface is
embedded, and since any two representatives of distinct closed curves are in minimal position, choosing $R$
sufficiently large, the above proof can be carried out verbatim.
\end{proof}

For convenience, we isolate the key idea involved here as it will be used again.

\begin{proposition} \label{lemma:main surgery argument}
If $f$ is a pseudo-Anosov element of $\Mod(S_g)$ with the property
that $\gin(c,f(c)) \geq n \geq 3$ for every simple closed curve $c$, then
$$\log(\lambda(f)) > \log\left(\frac{n}{2}\right).$$
\end{proposition}

\begin{proof}
As in the proof above, fix the metric $q = q_f$ on $S$, and let $c$ be a shortest curve in $S$ with respect to $q$. We
again assume geodesics are embedded and pairs are in minimal position, with the general case handled as above.
Take the shortest segment $a$ of $f(c)$ cut by $c$ (which is one of $\gin(c,f(c))$ segments of $f(c)$),
and the shortest segment $b$ of $c$ cut by $a$, and we obtain
$$\ell_q(c) \leq \ell_q(a \cup b) \leq \frac{\ell_q(f(c))}{\gin(c,f(c))} + \frac{\ell_q(c)}{2} < \frac{\lambda(f)
\ell_q(c)}{\gin(c,f(c))} + \frac{\ell_q(c)}{2}.$$ Dividing the left and right by $\ell_q(c)$, and simplifying and taking
logarithms, we obtain
$$\log(\lambda(f)) > \log \left( \frac{\gin(c,f(c))}{2} \right) \geq \log\left(\frac{n}{2}\right).$$
\end{proof}

\remark Wolpert \cite{Wo} has shown that a $K$-quasiconformal map $f$ of $S$ with respect to a hyperbolic metric $X$
distorts lengths in $X$ by a factor of at most $K$.  That is,
$\ell_X(f(c))/\ell_X(c) < K$, where $\ell_X(c)$ is the length of $c$
with respect to $X$.
In a previous
version we used this result and the same argument above (with no need for the final comment on minimal position and approximating Riemannian metrics) to produce a
lower bound of $.197$ for $\log(\lambda^2)$.
J. Franks suggested using
the quadratic differential metric, thus improving the lower bound by a
factor of $2$.

\subsection{Examples with small dilatation}
\label{section:upper}

In this section we give an upper bound for $L(\I(S))$ by constructing, for every $S=S_g$ ($g \geq 2$), an element $f
\in \I(S)$ with $\log(\lambda(f))< 4.127$.  We do this by appealing to a general construction for pseudo-Anosov mapping
classes given by Thurston \cite[\S6]{Th}; we refer the reader to that paper for the notation and details of the
construction.

A \emph{multicurve} is the isotopy class of a collection of pairwise disjoint simple closed curves, and a
\emph{multitwist} is the product of Dehn twists about the curves in a multicurve.

We begin by fixing a pair of multicurves $A = a_1 \cup \cdots \cup a_{\lceil g/2 \rceil}$ and $B = b_1 \cup \cdots
\cup b_{\lceil g/2 \rceil}$ in $S$ with the following three properties:
\begin{enumerate}
\item $A\cup B$ fills $S$.

\item $\gin(a_i,b_i) = \gin(a_i,b_{i-1}) = 4$ and
$\gin(a_i,b_j) = 0$ if $|i-j|>2$ (indices taken modulo $\lceil g/2 \rceil$).

\item Each $a_i$ and $b_j$ is a separating curve.
\end{enumerate}

\begin{figure}[htb]
\begin{center}
\psfrag{a1}{$\bar a_1$}
\psfrag{a2}{$\bar a_2$}
\psfrag{b1}{$\bar b_1$}
\psfrag{b2}{$\bar b_2$}
\includegraphics[scale=1]{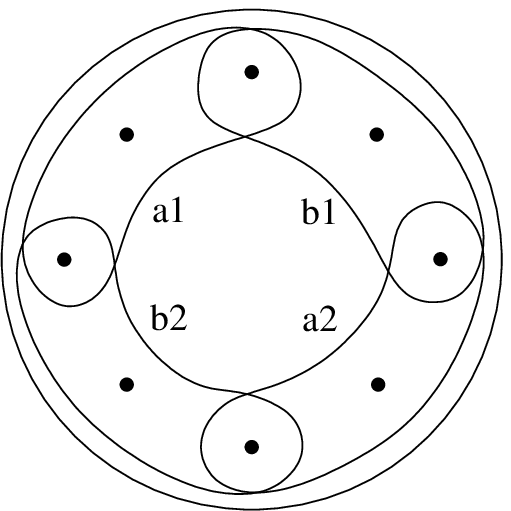} \hspace*{2cm} \includegraphics[scale=1]{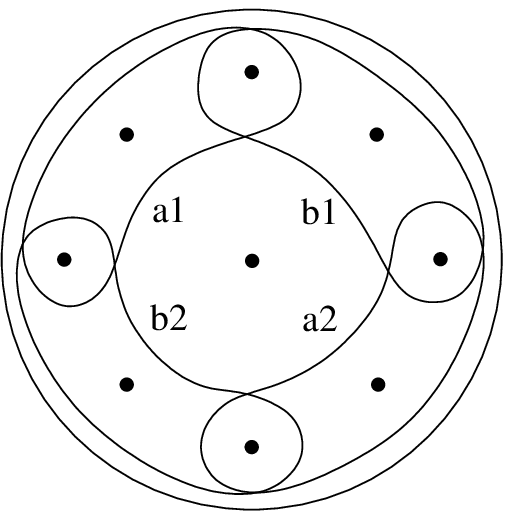}
\caption{$\bar{A}$ and $\bar{B}$ for $g = 3$ (left) and $g = 4$ (right; one marked point is ``in back'').}
\label{figure:small torelli}
\end{center}
\end{figure}

We can construct such an $A$ and $B$ explicitly as follows. Start with a sphere with $2g + 2$ marked points arranged
symmetrically as in Figure~\ref{figure:small torelli}; the arrangement depends on whether $g$ is odd (on the left) or
even (on the right---there is one more marked point ``in back''). Let $\bar{A}= \cup \bar a_i$ and $\bar{B} = \cup \bar
b_i$ be multicurves in the marked sphere as shown, and let $S$ be the two-fold cover, branched over the marked points,
with $A$ and $B$ the preimages of $\bar{A}$ and $\bar{B}$, respectively.  Since each component of $\bar{A}$ and
$\bar{B}$ surrounds exactly three marked points, each component of $A$ and $B$ is separating; in fact it bounds a genus
$1$ subsurface.

Next, we consider the matrix $N_{ij} = \gin(a_i,b_j)$, and compute the matrix $NN^t$.
This has entries given by
\[ (NN^t)_{ij} = \sum_{k=1}^{\lceil g/2 \rceil} \gin(a_i,b_k)\gin(a_j,b_k) \]
and the above description of intersection numbers easily implies that for $i$ and $j$ modulo $\lceil g/2 \rceil$ we
have
\[ (NN^t)_{ij} = \left\{ \begin{array}{rl}
32 &\quad \mbox{ for } \quad i = j\\
16 &\quad \mbox{ for } \quad |i - j| = 1\\
0 &\quad \mbox{ for } \quad |i - j| \geq 2.\\ \end{array} \right. \] In particular, note that the row sum of any row is $64$. It follows
that the Perron--Frobenius eigenvalue is $64$: take as an eigenvector the vector with all entries equal to $1$.

Now let $T_A$ denote the multitwist which is the composition of the
Dehn twists about each of the $a_i$ and $T_B$ the composition of Dehn
twists about each of the $b_j$.  In Thurston's construction of pseudo-Anosov homeomorphisms mentioned above
he begins by defining a homomorphism
$\langle T_A,T_B \rangle \to \PSL_2(\R)$ which in this case is given by
\[ T_A \mapsto \left( \begin{array}{rr} 1 & 8\\ 0 & 1\\
\end{array}
\right) \quad
\mbox{ and } \quad T_B \mapsto \left( \begin{array}{rr} 1 & 0\\ -8 &
1\\ \end{array} \right).\]
He then proves that any element of $\langle T_A,T_B \rangle$ that maps to a
hyperbolic element of $\PSL_2(\R)$ is pseudo-Anosov.
Moreover, the dilatation of such an element
is given by the absolute value of the leading eigenvalue of its image.

\remark
Thurston's construction of pseudo-Anosov homeomorphisms is much more general, merely requiring $A$ and $B$ to fill $S$.  In the general construction, the nonzero off-diagonal
entries of the homomorphic images of $T_A$ and $T_B$ are given by the square root of the Perron--Frobenius eigenvalue of $NN^t$.  Again, $N$ is the matrix of intersection numbers of components of $A$ and $B$.\\

In our case, a direct computation shows that the mapping class $f = T_AT_B$ maps to a matrix with trace $=-62$.  It
follows that $f$ is pseudo-Anosov, and that $\lambda = \lambda(f)$ satisfies
\[ \lambda^2 - 62 \lambda + 1=0. \]
Solving for the largest root, we find $\log(\lambda) < 4.127$.  Since Dehn twists about separating curves are elements of $\I(S)$, we have proven the upper bound of
Theorem~\ref{theorem:torelli}.

\begin{proposition}
\label{proposition:torelli upper bound} If $g \geq 2$, then $L(\I(S_g)) < 4.127$.
\end{proposition}

\begin{remark}
The map $T_A T_B$ has the smallest dilatation among all pseudo-Anosov
elements of $\langle T_A, T_B \rangle$ (see \cite{Le}).
\end{remark}

\medskip
This construction also provides a universal upper bound for $L(PB_n)$
for $n\geq 3$.

\begin{theorem}
For all $n \geq 3$, we have
$$1.443< L(PB_n) < 2.634.$$
\end{theorem}

\begin{proof}
As mentioned in the introduction, Song \cite{So} proved the lower
bound $1.443 \approx \log(2+ \sqrt{5})$.  For the upper bound, we start with the case where $n$ is odd, say $n=2g+1$.  To prove the upper bound, consider the
sphere with marked points which we described above.  We can puncture
every marked point and turn one puncture into a boundary component,
making the surface into a $(2g+1)$-times punctured disk.  Then
$T_{\bar A}T_{\bar B}$ represents a pseudo-Anosov braid in
$PB_{2g+1}$.  Indeed, we can use the same method of Thurston described
above to find the homomorphism $\langle T_{\bar A},T_{\bar B} \rangle
\to \PSL_2(\R)$.  The eigenvalue for the analogous $NN^t$ matrix is
16, and so by a calculation, we obtain the upper bound
\[
\log( \lambda(T_{\bar{A}}T_{\bar{B}})) < 2.634.
\]

We now do the case when $n$ is even.  Whenever a marked point is not
contained in a bigon, we can erase the marking, and $\bar A$ and $\bar
B$, as drawn, are still in minimal position.  Puncturing the remaining
marked points and turning one puncture into a boundary component, we can
obtain examples proving the upper bound $L(PB_n) < 2.634$ for $n
\geq 3$.
\end{proof}

\subsection{Principal congruence subgroups}
\label{section:principal congruence}

We now give the proof of Theorem~\ref{theorem:principal congruence}, which states that the bounds given in
Propositions~\ref{proposition:torelli lower bound} and~\ref{proposition:torelli upper bound} for $L(\I(S))$ can be
extended to $\Mod(S)[r]$ when $r \geq 3$.

\begin{question}
Is it true that $L(\Mod(S_g)[2]) \asymp 1$?
\end{question}

\begin{proof}[Proof of Theorem~\ref{theorem:principal congruence}]

Since $\I(S) < \Mod(S)[r]$, the upper bound is immediate.  For $r \geq 4$, the proof of the lower bound is essentially
the same as that of Proposition~\ref{proposition:torelli lower bound}; all that needs to be verified is that
Lemmas~\ref{lemma:separating curves intersect image} and~\ref{lemma:nonseparating curves intersect image} hold under
the weaker hypothesis that $f \in \Mod(S)[r]$, $r \geq 4$. Indeed, the same arguments work, using splittings of
$H_1(S;\Z/r\Z)$ (and its quotients) in place of $H_1(S;\Z)$.

For the case $r=3$, the proof is the same as in the $r \geq 4$ case, except we need to include the possibilities
$\gin(c,f(c)) = 3$ and $\gin(c,f^2(c)) = 3$ in Case 1 of Proposition~\ref{proposition:torelli lower bound}.  By
Proposition~\ref{lemma:main surgery argument}, the lower bound for this case becomes $.202
> .197$, and the argument for Case 2 still gives a lower bound of $.197$, so we are done.
\end{proof}

It follows from the discreteness of $\spec(\Mod(S))$, and the fact that $\I(S) < \Mod(S)[r]$, that there is a (minimal)
$r=r(g)$ such that $L(\Mod(S_g)[n]) = L(\I(S_g))$ whenever $n \geq r$ (see the proof of
Proposition~\ref{proposition:filtration}).

\begin{question}
What are the values of $r(g)$?  What are the asymptotics of $r(g)$?
\end{question}

\section{The Johnson kernel}
\label{section:johnson kernel}

Johnson \cite{Jo1} proved that $\K(S_g)$ is an infinite index subgroup of $\I(S_g)$ for $g \geq 3$ (when $g=2$, the
two groups agree). We have $L(\K(S))\geq L(\I(S))$ since $\K(S) < \I(S)$, and it is natural to ask the following.

\begin{question}
Is $L(\K(S_g)) > L(\I(S_g))$ for $g \geq 3$?
\end{question}

While we do not know the answer to this question, we are able to give a
better lower bound for $L(\K(S))$ than we did for $L(\I(S))$ in
Theorem~\ref{theorem:torelli}.  As with $\I(S)$, the key is to understand
how elements of $\K(S)$ act on curves.

\subsection{\boldmath$\K(S)$ and geometric intersection numbers}

The conclusions of Lemmas~\ref{lemma:separating curves intersect image} and~\ref{lemma:nonseparating curves intersect
image} can be improved by assuming $f \in \K(S)$.

\begin{proposition}
\label{lemma:johnson intersections}
For $f \in \K(S)$, and any curve $c$, if $c \neq f(c)$, then
$\gin(c,f(c)) \geq 4$.
\end{proposition}

The proposition is sharp, since for $g \geq 2$ one can find a curve
$c$ and a separating curve $d$ with $\gin(c,d) = 2$, and in this case
$\gin(T_d(c),c) = 4$.

\begin{proof}
When $c$ is separating the proposition was already proven in Lemma~\ref{lemma:separating curves intersect image} for any
$f\in \I(S)$.  So assume that $c$ is nonseparating. Since $f(c)$ is homologous to $c$, it suffices to rule out
$\gin(c,f(c))=0$ and $\gin(c,f(c))=2$. As $\K(S)$ is normal in $\Mod(S)$, the mapping class
\[ [T_c,f] = T_cfT_c^{-1}f^{-1}=T_cT_{f(c)}^{-1} \]
must lie in $\K(S)$.  The proposition now follows from Lemma~\ref{lemma:no xmaps in johnson} below.
\end{proof}

\begin{lemma}
\label{lemma:no xmaps in johnson}
$\K(S)$ contains no elements of the form $T_cT_d^{-1}$ where $c$ and $d$ are distinct homologous curves with $\gin(c,d)$ either 0 or 2.
\end{lemma}

\begin{proof}

Johnson \cite{Jo1} constructed a homomorphism $\tau:\I(S)\to
(\wedge^3H)/(\langle \omega \rangle \wedge H)$, where $H=H_1(S;\Z)$ and
$\omega$ is the symplectic intersection pairing.  He proved that the
kernel is exactly $\K(S)$.  Moreover, Johnson \cite[Corollary to Lemma
4B]{Jo1} gave an explicit formula for the $\tau$-image of a {\em bounding
pair map}, i.e.\ a product of twists $T_aT_b^{-1}$ where $\{a,b\}$ is a
bounding pair.  For the formula, let $R$ be the component of $S-(a \cup b)$ not containing the base point for $\pi_1(S)$, let $[a]$ denote the homology class of $a$ and $b$ (oriented so $R$ is to the left of $a$), and let $u_1,v_1
\dots, u_k, v_k$ is any symplectic basis for $H_1(R;\Z)/\langle [a]
\rangle$.  The formula reads:
\begin{equation}
\label{equation:tau of bp}
 \tau(T_aT_b^{-1}) = \left(\sum_{i=1}^k u_i \wedge v_i\right)\wedge [a].
\end{equation}
It immediately follows that $\K(S)$ contains no bounding pair maps, and
hence it remains to show that $\K(S)$ contains no elements of the form
$T_cT_d^{-1}$ where $c$ and $d$ are homologous curves with
$\gin(c,d)=2$.  By Lemma~\ref{lemma:homologous curves git2}, $c$ and $d$ are necessarily configured as
in Figure~\ref{figure:xpairs}.  Using the notation of the picture, the
lantern relation (see \cite[\S7g]{De}) gives
\[ T_e T_c T_d  = T_x T_y T_z T_w \]
which implies
\[ T_c T_d^{-1} = T_e^{-1} T_x T_y T_z T_w T_d^{-2}. \]
As $T_e$, $T_x$, and $T_y$ are elements of $\K(S)$, we see that $T_cT_d^{-1}$ is an element of $\K(S)$ if and only if $T_zT_wT_d^{-2}$ is an element of $\K(S)$.  The latter is a product of two bounding pair maps: $(T_zT_d^{-1})(T_wT_d^{-1})$.  To prove the lemma then, it suffices to check that $\tau(T_zT_d^{-1}) \neq \tau(T_dT_w^{-1})$.
But this is apparent from equation~(\ref{equation:tau of bp}) (consult Figure~\ref{figure:xpairs}).
\begin{figure}[ht]
\begin{center}
\psfrag{x}{$x$} \psfrag{y}{$y$} \psfrag{z}{$z$} \psfrag{w}{$w$} \psfrag{c}{$c$} \psfrag{d}{$d$} \psfrag{e}{$e$}
\includegraphics[scale=.9]{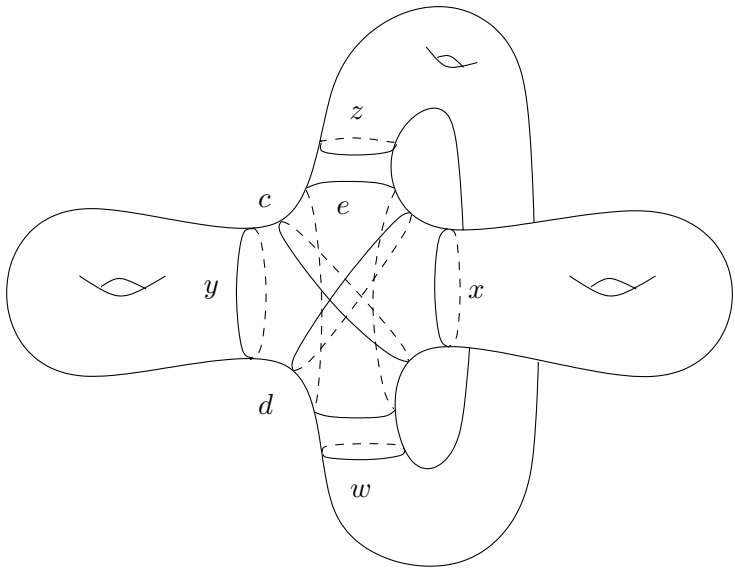}
\caption{Homologous curves $c$ and $d$ with $\gin(c,d) = 2$.  The genera of the subsurfaces bounded by $x$, $y$, and
the pair $\{z,w\}$ may vary.} \label{figure:xpairs}
\end{center}
\end{figure}
\end{proof}

\subsection{Bounds for $L(\K(S))$}

We are now ready to give the following improvement for the
bounds on $L(\K(S))$ given by Theorem~\ref{theorem:torelli}.

\begin{proposition}
\label{proposition:johnson} For $g \geq 2$, we have $$.693< L(\K(S_g))< 4.127.$$
\end{proposition}

\begin{proof} The mapping class $T_AT_B$ constructed in \S\ref{section:upper}
is a composition of Dehn twists about separating curves.  Thus this mapping class already lies in $\K(S)$, giving the
upper bound.  The lower bound follows immediately from Propositions~\ref{lemma:johnson intersections}
and~\ref{lemma:main surgery argument}, as $\log(2) \approx .693$ approximated from below.
\end{proof}

\section{The Johnson filtration}
\label{section:filtrations}

We will now prove Theorem \ref{theorem:johnson filtration}.  Both the
upper and lower bounds will follow from generalized versions of our
arguments for $\I(S)$ and $\K(S)$.

\subsection{Asymptotic lower bounds}

Before proving the lower bound in Theorem \ref{theorem:johnson filtration}, we give the following weaker statement
which holds for any \emph{normal filtration} of $\Mod(S)$, by which we mean a filtration of $\Mod(S)$ by normal
subgroups.

\begin{proposition}
\label{proposition:filtration}
For any normal filtration $N_1 > N_2 > \cdots$ of $\Mod(S)$, we have $L(N_k) \to \infty$ as $k \to \infty$.
\end{proposition}

\begin{proof}
Given $M > 0$, there are only finitely many
conjugacy classes of pseudo-Anosov mapping classes $f$ with $\lambda(f) \leq M$
(see \cite{Iv1}).
The proposition then follows from the definition of a normal filtration.
\end{proof}

The first step towards the proof of Theorem~\ref{theorem:johnson
filtration} is to generalize Lemma~\ref{lemma:no xmaps in johnson}.
In the proof of this lemma, it was essential that there were unique
pictures for homologous curves with geometric intersection number 0 or 2.  We are forced to
replace this precise description of how our two curves sit in $S$ with a
rough finiteness statement.

A \emph{configuration} is a triple $(S,c,d)$, where $S$ is a closed
surface, and $c$ and $d$ are distinct curves in $S$ which are in minimal
position (i.e. their union does not bound any bigon). Let $N =
N(c,d)$ denote a closed regular neighborhood of $c \cup d$. There is a
natural partial ordering on configurations where $(\hat S, \hat c,\hat
d) < (S,c,d)$ if $\hat S \not \cong S$ and there is a continuous map
$\eta:(S,c,d) \to (\hat S,\hat c,\hat d)$ which restricts to a
homeomorphism of triples $$\eta|_{N(c,d)} : (N(c,d),c,d) \to (N(\hat
c,\hat d),\hat c, \hat d).$$ We call
such a map $\eta$ a \emph{crushing map}.  Because the composition of
crushing maps is a crushing map, $<$ is a partial order.  Any
minimal configuration with respect to this partial ordering is called a
\emph{terminal configuration}. We declare two configurations $(S,c,d)$
and $(\hat S, \hat c, \hat d)$ to be equal provided they are
homeomorphic as triples.

\begin{lemma}
\label{lemma:finitely many configurations}
For every $n> 0$ there are only finitely many terminal configurations $(S,c,d)$ with $\gin(c,d) \leq n$.
\end{lemma}

\begin{proof}  To begin, we show that the topology of the complement of $N(c,d)$ for any terminal configuration $(S,c,d)$ is limited.\\

\noindent Claim: If $(S,c,d)$ is terminal then every complementary component $U$ of $S - N(c,d)$ has genus at most one.\\

\noindent Proof of claim: Suppose $(S,c,d)$ is a configuration, and some component $U$ has genus at least two. We
produce a crushing map $\eta:(S,c,d) \to (\hat S,\hat c,\hat d)$ as follows. Let $\eta: S \to \hat S$ be the quotient
of $S$ obtained by first fixing a compact genus 1 subsurface $R \subset U \subset S$ with exactly one boundary
component and identifying (``crushing'') $R$ to a point.  We let $\hat c = \eta(c)$ and $\hat d = \eta(d)$ and note
that the restriction of $\eta$ to $N(c,d)$ is a homeomorphism onto $N(\hat c,\hat d)$.  To prove that $\eta$ is a
crushing map, all that remains is to verify that $(\hat S,\hat c,\hat d)$ is indeed a configuration.  That is, we must
check that $\hat c$ and $\hat d$ are essential, not homotopic to one another, and in minimal position.  For this, it
suffices to verify that no component of the complement of $N(\hat c,\hat d)$ is a disk or annulus if $\hat c \cap \hat
d = \emptyset$ or a bigon if $\hat c \cap \hat d \neq \emptyset$.  However, the components of the complement of $N(\hat
c,\hat d)$ are all homeomorphic to those of $N(c,d)$ with the exception of $\eta(U)$, and since $c$ and $d$ are
essential, homotopically distinct, and in minimal position, it suffices to verify this statement for the single
component $\eta(U)$. By construction, $\eta(U)$ has genus at least one, so it is not a disk, annulus, or bigon, and
therefore $(\hat S, \hat c, \hat d)$ is a configuration.
It follows that $(\hat S,\hat c,\hat d) < (S,c,d)$, and $(S,c,d)$ is not terminal, proving the claim.\\

Now let $(S,c,d)$ be a terminal configuration with $N = N(c,d)$.
By Lemma \ref{lemma:filling curves intersect}, we have
$\chi(N) = - \gin(c,d) \geq -n$, and so there are finitely many possibilities for $N$, up to homeomorphism. Since $c$
is a curve in $N$, there are only finitely many possibilities for the homeomorphism type of $(N,c)$. Further, since $c$
cuts $d$ into $\gin(c,d) \leq n$ arcs, it follows that there are only finitely many possibilities for the homeomorphism
type of $(N,c,d)$.

Now note that there are only finitely many possibilities for the number of boundary components of $N$, and hence
finitely many possibilities for the number of boundary components of $\overline{S - N}$.  Because each component of
$\overline{S - N}$ has genus at most $1$, there are only finitely many possibilities for the homeomorphism type of
$\overline{S - N}$.  Finally, the homeomorphism type of $(S,c,d)$ can be specified by $\overline{S - N}$ and $(N,c,d)$
and the (finite) combinatorial gluing data matching boundary components of the former with those of the latter.
\end{proof}

The next lemma allows us to say that, given a particular terminal configuration $(S,c,d)$ with $T_cT_d^{-1} \in \N_k(S)$, it is not possible to
``push'' $T_cT_d^{-1}$ further down the Johnson filtration by adding genus to $S$ outside $N(c,d)$.  The proof applies
in a much more general context, so we state it in this generality.

Suppose we are given a map of pairs $\eta:(S,N) \to (\hat S,\hat N)$ where $N \subset S$ and $\hat N \subset \hat S$
are subsurfaces and the restriction $\eta|_N : N \to \hat N$ is a homeomorphism; for example $\eta$ might be a crushing
map.  Then any $f \in \Mod(S)$ which is supported in $N$ pushes forward via $\eta$ to an element $\hat f \in \Mod(\hat
S)$ supported in $\hat N$ given by $f|_{\hat N} = \eta|_N \circ f|_N \circ \eta|_N^{-1}$.

\begin{lemma}
\label{lemma:handle crushing} Suppose $\eta:(S,N) \to (\hat S,\hat N)$, $f \in \Mod(S)$, and $\hat f \in \Mod(\hat S)$
are as above. Then $\hat f \in \N_k(\hat S)$ whenever $f \in \N_k(S)$.
\end{lemma}

\begin{proof}

Suppose that $f \in \N_k(S)$, i.e. after picking a representative of $f$ and fixing a base point, the induced action
$f_\star$ of $f$ on $\Gamma/\Gamma_k$ is inner (where $\Gamma = \pi_1(S)$ and $\Gamma_k$ is the $k^{\mbox{\tiny th}}$
term of its lower central series, as above).

Let $\hat \Gamma=\pi_1(\hat S)$, and denote by $\{\hat \Gamma_i\}$ its lower central series.  The map $\eta$ induces a
surjective homomorphism $\Gamma \to \hat \Gamma$ which restricts to a surjection $\Gamma_k \to \hat \Gamma_k$, so we
have an induced map $\eta_\star: \Gamma/\Gamma_k \to \hat \Gamma/\hat \Gamma_k$.  Finally, let $\hat f_\star$ be the
induced action of $\hat f$ on $\hat \Gamma/\hat \Gamma_k$.  We encode this information in the following diagram.
\begin{figure}[ht]
\begin{center}
\scalebox{1}{
\xymatrix{
\Gamma/\Gamma_k \ar@{->}[r]^{f_\star} \ar@{->}[d]_{\eta_\star} & \Gamma/\Gamma_k \ar@{->}[d]^{\eta_\star} \\
\hat \Gamma/\hat \Gamma_k \ar@{->}[r]^{\hat f_\star} & \hat \Gamma/\hat \Gamma_k \\
}
}
\end{center}
\end{figure}

The diagram is commutative by the definition of $\hat f_\star$, which implies that $\hat f_\star$ is also inner; indeed, if
$f_\star$ is conjugation by $\gamma$, then $\hat f_\star$ is conjugation by $\eta_\star(\gamma)$. Therefore, $\hat f
\in \N_k(\hat S)$.
\end{proof}

We now arrive at the desired generalization of Lemma~\ref{lemma:no xmaps in
johnson}.  Define
\[ C(n) = 1+\sup \{k \, | \, (S,c,d) \mbox{ is terminal, } \gin(c,d) \leq n \mbox{, and } T_c T_d^{-1} \in \N_k(S) \}
\]
By Lemma~\ref{lemma:finitely many configurations} and the definition of a normal filtration, $C(n)$ is finite for each
$n$.

\begin{lemma}
\label{lemma:deep implies high intersection} Let $n > 0$.  If
$c \neq d$ and $\gin(c,d)
\leq n$, then $T_cT_d^{-1} \notin \N_{C(n)}(S)$.  Furthermore,
$\displaystyle \lim_{n \to \infty} C(n) = \infty$.
\end{lemma}

\begin{proof}
Suppose $c$ and $d$ are curves in $S$ with $\gin(c,d) \leq n$ and $T_cT_d^{-1} \in \N_k(S)$.   If $(S,c,d)$ is a
terminal configuration, then by the definition of $C(n)$, we see that $k < C(n)$.  If $(S,c,d)$ is not terminal, then
there is a terminal configuration $(\hat S, \hat c, \hat d)$ and crushing map $\eta:(S,c,d) \to (\hat S, \hat c, \hat
d)$.  The induced map of pairs $\eta:(S,N(c,d)) \to (\hat S, N(\hat c, \hat d))$ and the map $T_cT_d^{-1}$ satisfy the
hypothesis of Lemma \ref{lemma:handle crushing} with
\[ \widehat{T_cT_d^{-1}} = T_{\hat c} T_{\hat d}^{-1} \]
so $T_{\hat c} T_{\hat d}^{-1} \in \N_k(\hat S)$, which again implies $k < C(n)$ as required.

To complete the proof and show that $C(n) \to \infty$, we first notice
that $C(n)$ is nondecreasing, since the set of terminal configurations
used to define $C(n)$ contains the set of configurations used to
define $C(n-1)$.  Thus, it suffices to show that for any $k$ there
exists a $c$ and $d$ such that $T_cT_d^{-1} \in \N_k(S)$ for some $S$.
Let $f \in \N_k(S)$ be any nontrivial element and $c$ any curve with
$f(c) \neq c$.  Since $\N_k(S) \lhd
\Mod(S)$, it follows that $[T_c,f]=T_cT_{f(c)}^{-1}$ is an element of
$\N_k(S)$.
\end{proof}

We are finally ready to give the ``asymptotic version'' of
Proposition~\ref{lemma:johnson intersections}.  For the statement,
define
\[ B(k) = \sup \{n:C(n) \leq k\} + 1. \]
In the case where $B$ is not defined by this equation, we artificially set $B = 0$ ($B$ is not defined for any integer which is smaller than the smallest value of $C$).
Note that $B(k)$ is well-defined and finite for each $k \geq 0$ since $C$ is unbounded and nondecreasing.
Rephrasing, $B(k)$ is the minimum intersection number required for any pair of curves $c$ and $d$ in any surface $S$ to satisfy $T_cT_d^{-1} \in \N_k(S)$.

We require the following alternate characterization of $B(k)$, which follows immediately from the definition.

\begin{lemma}
\label{lemma:B characterization}
$B(k)$ is the smallest integer valued function for
which $C(B(k)) > k$ for all $k$.
\end{lemma}

The next proposition
gives the desired generalization of Proposition \ref{lemma:johnson
intersections}.

\begin{proposition}
\label{lemma:asymptotic intersections} Let $k > 0$ and let $S$ be any surface.  If $f \in \N_k(S)$, then
$\gin(c,f(c)) \geq B(k)$ for every simple closed curve $c$ in $S$
with $f(c) \neq c$.  Moreover, $\displaystyle \lim_{k \to
\infty}B(k) = \infty$.
\end{proposition}

\begin{proof}
First, $B(k) \to \infty$ as $k \to \infty$ since $C$ is unbounded and
nondecreasing.  Now, given $k$, choose any $S$, any $f \in \N_k(S)$,
and any simple closed curve $c$ in $S$ with $f(c) \neq c$.  Since
$\N_k(S)$ is normal in $\Mod(S)$, we have $[T_c,f]=T_cT_{f(c)}^{-1}
\in \N_k(S)$.  By Lemma~\ref{lemma:B characterization}, if
$\gin(c,f(c)) < B(k)$, then $C(\gin(c,f(c))) \leq k$.  By
Lemma~\ref{lemma:deep implies high intersection}, we have
$T_cT_{f(c)}^{-1} \notin \N_k(S)$, which is a
contradiction.
\end{proof}

Using Propositions~\ref{lemma:asymptotic intersections}
and~\ref{lemma:main surgery argument}, it is now straightforward to
prove the lower bound of Theorem~\ref{theorem:johnson filtration}.

\begin{proof}[Proof of Theorem~\ref{theorem:johnson filtration}.]

If $f \in \N_k(S)$ is pseudo-Anosov, then $\gin(c,f(c)) \geq B(k)$ for every curve $c$, by Proposition~\ref{lemma:asymptotic
intersections} and the fact that a pseudo-Anosov mapping class does
not fix any curve.  We set
$$m(k) =
\log \left( \frac{B(k)}{2} \right).$$ By Proposition~\ref{lemma:main surgery argument}, $\log(\lambda(f)) > m(k)$, and
so this completes the proof.
\end{proof}

Among several questions which now arise, we pose the following.

\begin{question}
What are the asymptotics of $B(k)$?  What are the asymptotics of $L(\N_k(S))$?
\end{question}

\subsection{Asymptotic upper bounds}
\label{section:asymptotic upper bounds}

We now prove the upper bound in Theorem \ref{theorem:johnson
filtration}.

\begin{proposition}
For any $k \geq 1$, there is an $M(k)$ so that $L(\N_k(S_g)) < M(k)$
for all $g \geq 2$.
\end{proposition}

Since $m(k) \to \infty$, it follows that $M(k) \to \infty$ as $k \to \infty$.

\begin{proof}

Let $k \geq 1$ be fixed.  To prove the proposition, we need to find a pseudo-Anosov mapping class $f\in\N_k(S)$ whose
dilatation depends on $k$, but not on $S$.  We begin by recalling that since $\{ \N_i(S) \}_{i \geq 1}$ is a central
series for $\I(S)$, then the $(k-1)^{\mbox{\tiny st}}$ term of the lower
central series of $\I(S)$ is contained in $\N_k(S)$; note
that $\N_1(S) = \I(S)$ is the zero term of the lower central series.

Now, without specifying a particular surface $S$, we consider the group $\langle T_A, T_B \rangle$ generated by the
multitwists $T_A$ and $T_B$ of Section~\ref{section:upper}.  The group $\langle T_A, T_B \rangle$ is a free group on
the given generators (see \cite[\S 6.1]{Le} for a discussion).  Therefore, there is a nontrivial element $f$ in the
$(k-1)^{\mbox{\tiny st}}$ term of the lower central series of $\langle T_A, T_B \rangle$.  Since $T_A$ and $T_B$ are
both elements of $\I(S)$, it follows that $f$ is an element of the $(k-1)^{\mbox{\tiny st}}$ term of the lower central
series of $\I(S)$, and hence $f \in \N_k(S)$.

The key feature here is this: the image of $f$ in $\PSL_2(\R)$ does not depend of the choice of $S$.  This is
because $f$ was chosen independently of $S$ as a word in $T_A$ and $T_B$, and the images of $T_A$ and $T_B$ in
$\PSL_2(\R)$ do not depend on the choice of $S$ (see Section~\ref{section:upper}).

Since $\langle T_A, T_B \rangle$ is a free group and the only elements of this group which are not pseudo-Anosov are
conjugates of $T_A$ and $T_B$ (see, e.g., \cite{Le}), it follows that $f$ is pseudo-Anosov.  Since its dilatation only
depends on its image in $\PSL_2(\R)$, and the latter is independent of the choice of $S$, we are done.
\end{proof}

\remark Note that the word in $T_A$ and $T_B$ given as a simple nested commutator has word length on the order of
$2^k$, where $k$ is the number of nested commutators involved
(i.e. the depth in the lower central series).  Thus the order of logarithm of the dilatation is at most exponential
in $k$.

\section{Translation lengths on the complex of curves}
\label{section:complex of curves}

Our goal in this section is to prove Theorems~\ref{theorem:curve complex infs} and~\ref{theorem:LCC curve
complex}.  These will follow rather quickly from Theorem~\ref{theorem:curve complex v. dilatation}.

We first need the following technical fact.

\begin{lemma}
\label{lemma:liminf}
If $m,n \in \Z$ and $f \in \Mod(S)$ satisfy $n\tau_\C(f) > m$, then $d_\C(f^n(c)) \geq m+1$ for any curve $c$.
\end{lemma}

\begin{proof}

If not, we have $d_\C(f^n(c),c) \leq m$, so by the triangle inequality
$d_\C(f^{nj}(c),c) \leq mj$. Dividing both sides by $nj$ and taking the $\liminf$, we get
\[ \liminf_{j \to \infty} \frac{d_\C(f^{nj}(c),c)}{nj} \leq \frac{m}{n}. \]
The $\liminf$ used to define $\tau_\C(f)$ is no larger than the left hand side, and so we arrive at $n\tau_\C(f) \leq
m$, a contradiction.
\end{proof}

\begin{theorem}
\label{theorem:curve complex v. dilatation} For any $g \geq 2$ and any pseudo-Anosov $f \in \Mod(S_g)$ with $\lambda(f)
\leq g - 1/2$, we have
$$\tau_\C(f) < \frac{4 \log(\lambda(f))}{\log\left(g-\frac{1}{2}\right)}.$$
\end{theorem}
\begin{remark} It seems likely that the hypothesis $\lambda(f) \leq g - 1/2$ is not necessary, but it is required for
our argument.
\end{remark}

\begin{proof}
Let $n$ be the smallest integer so that $2 < n \tau_\C(f)$.  Note that $n\tau_\C(f) \leq 4$ whenever $n > 1$.

Now, let $c$ be any curve in $S$.  By Lemma~\ref{lemma:liminf}, $d_\C(f^n(c),c) \geq 3$, which (by the definition of $d_\C$)
implies that $c$ and $f^n(c)$ fill $S$.  Lemma \ref{lemma:filling curves intersect} implies $\gin(c,f^n(c)) \geq 2g - 1$, and hence Proposition~\ref{lemma:main
surgery argument} applied to $f^n$ says
\[ n \log(\lambda(f)) = \log(\lambda(f^n)) > \log \left( g - \frac{1}{2} \right) \]
which we write as
\[ \frac{1}{n} < \frac{\log(\lambda(f))}{\log\left(g - \frac{1}{2}\right)}. \]
By hypothesis, the right hand side is at most 1, and so $n > 1$.  As mentioned above, this means that $n\tau_\C(f) \leq
4$. Thus, we have
\[ \tau_\C(f) \leq \frac{4}{n} < \frac{4 \log(\lambda(f))}{\log\left(g - \frac{1}{2}\right)}. \]
\end{proof}

We can now deduce Theorems~\ref{theorem:curve complex infs} and~\ref{theorem:LCC curve complex} as corollaries of
Theorem~\ref{theorem:curve complex v. dilatation}.

\begin{proof}[Proof of Theorem~\ref{theorem:curve complex infs}]
Let $f_g \in \Mod(S_g)$ be a minimal dilatation pseudo-Anosov mapping class.  Hironaka--Kin \cite{HK} showed that
$\log(\lambda(f_g)) \leq \log(2+ \sqrt 3)/g$, and so if $g \geq 3$, then $\lambda(f_g) < g- 1/2$. The theorem thus
follows for $g \geq 3$ from Theorem \ref{theorem:curve complex v. dilatation}.  The case of genus 2 can be handled by
explicit examples.
\end{proof}

\begin{proof}[Proof of Theorem~\ref{theorem:LCC curve complex}]
For any fixed $k$, with $M(k)$ as in Theorem \ref{theorem:johnson filtration}, we have $M(k) \leq \log(g- 1/2)$ for $g$
sufficiently large.  That is, for large enough $g$, we have some $f_g \in \N_k(S_g)$ with $\lambda(f_g) \leq g - 1/2$.
Letting $g$ tend to infinity, Theorem~\ref{theorem:curve complex v. dilatation} implies
\[ L_\C(\N_k(S_g)) \leq \tau_\C(f_g) < \frac{4 \log(\lambda(f_g))}{\log \left( g - \frac{1}{2} \right) } \leq \frac{4 M(k)}{\log \left( g -
\frac{1}{2} \right) } \to 0. \]
\end{proof}

\section{Brunnian subgroups}
\label{section:brunnian}

Let $S_{g,p}$ be the orientable surface of genus $g$ with $p>0$ punctures, and let $\PMod(S_{g,p})$ be the subgroup of
$\Mod(S_{g,p})$ consisting of elements which fix each puncture. There are $p$ natural surjective homomorphisms
\[ F_i: \PMod(S_{g,p}) \to \PMod(S_{g,p-1}) \]
obtained by filling in the $i^{\mbox{\tiny th}}$ puncture, for $1 \leq i
\leq p$.  The \emph{Brunnian subgroup} of
$\Mod(S_{g,p})$ is the (nonempty!) intersection of the kernels:
\[ \brun(S_{g,p}) = \displaystyle \bigcap_{i=1}^p \textrm{ker}(F_i). \]
A topological description of each $F_i$ is given by the Birman exact
sequence \cite[Theorem 1.4]{Bi}.

\begin{proof}[Proof of Theorem \ref{theorem:brunnian}] This is similar to
the proof of the lower bound in Proposition \ref{proposition:torelli lower bound} and comes in two parts. We begin by
uniformly bounding $\gin(c,f(c))$ from below for any $f \in \brun(S_{g,p})$ and any curve $c$ with $f(c) \neq c$. To do
this, we first note that by definition $F_i(f)(c)=c$ for every curve $c$ and every $i=1,...,p$ (since $F_i(f)=1$). In
other words, if we fill in any puncture, $f(c)$ becomes isotopic to $c$.  Therefore, the complement of $c \cup f(c)$
contains $p$ punctured bigons, one for each puncture of $S_{g,p}$.  In the present case, an endpoint of a punctured bigon
can lie in at most two punctured bigons and so we have $\gin(c,f(c)) \geq p$.

For the second part of the proof we would like to apply Proposition~\ref{lemma:main surgery argument}.  However, this
is unavailable: the hypothesis of that proposition requires that the surface involved be closed.  Indeed, that proof
breaks down when the surface has punctures since the curve which is produced by the cut-and-paste may be
\emph{peripheral} (homotopic to a puncture), and hence has no geodesic representative.

Proposition~\ref{lemma:punctured surgery argument} below is a version of Proposition~\ref{lemma:main surgery argument}
for punctured surfaces, and it completes the proof.
\end{proof}

\begin{proposition} \label{lemma:punctured surgery argument}
If $f \in \Mod(S_{g,p})$ is pseudo-Anosov and has the property that
$\gin(c,f(c)) \geq n \geq 5$ for every simple closed curve $c$, then
$$\log(\lambda(f)) > \log\left(\frac{n}{4}\right).$$
\end{proposition}

\begin{proof}
As in the proof of Proposition~\ref{lemma:main surgery argument}, we let
$q=q_f$.  In addition to the fact that the metric is singular and so
geodesic representatives may not be in minimal position, there is
another difficulty which arises in this setting.  Namely, the presence
of punctures makes the metric incomplete and geodesic representatives
may not exist at all.  We modify the metric to be a complete Riemannian
metric to alleviate both of these problems.  As in the proof of
Proposition~\ref{proposition:torelli lower bound}, we can change the
metric in small neighborhoods of the singularities to be smooth and have
nonpositive curvature.  We can also modify the metric in a small
neighborhood of the punctures to be nonpositively curved and complete by
inserting a hyperbolic cusp and interpolating between the hyperbolic
metric and flat metric by nonpositively curved metrics (again, by
explicit computation).

Let $q_0$ denote the modified Riemannian metric of nonpositive curvature.  We may thus assume that all $q_0$--geodesics are embedded and pairs are in
minimal position.  Moreover, by choosing $q_0$ to approximate $q$ sufficiently well on large compact subsets of $S$, we may assume that for all sufficiently short nonperipheral curves (in particular, all those curves that we will encounter)
\[
\frac{\ell_{q_0}(f(c))}{\ell_{q_0}(c)} < \lambda(f).
\]

We let $c$ be a shortest (nonperipheral) curve in the $q_0$-metric and consider two arcs $a_1$ and $a_2$ of $f(c)$ cut
along $c$ which share an endpoint, and for which
\[
\ell_{q_0}(a_1) + \ell_{q_0}(a_2) \leq 2~\frac{\ell_{q_0}(f(c))}{\gin(c,f(c))}.
\]
Let $b_1$ and $b_2$ be the shortest arcs of $c$ cut by $a_1$ and $a_2$, respectively.  We also consider the
concatenated arc $a = a_1 \cup a_2$, and let $b$ denote the shortest arc of $c$ cut by $a$.

Suppose now that $a_1 \cup b_1$, say, is not peripheral. Then as in the proof of Proposition~\ref{lemma:main surgery
argument} we obtain
\[
\ell_{q_0}(c) \leq \ell_{q_0}(a_1 \cup b_1) < \frac{2\lambda(f) \ell_{q_0}(c)}{n} + \frac{\ell_{q_0}(c)}{2}
\]
and hence
\[
\log(\lambda(f)) > \log\left(\frac{n}{4}\right).
\]
Since each of $a_1$, $a_2$, and $a$ has length at most $2\ell_{q_0}(f(c))/\gin(c,f(c))$, and each of $b_1$, $b_2$, and
$b$ has length at most $\ell_{q_0}(c)/2$, we obtain the same bound if any of $a \cup b$, $a_1 \cup b_1$, or $a_2 \cup
b_2$ is nonperipheral. Thus the proof will be complete if we can show that this is the case.

We label the endpoints of $a_1$ and $a_2$ as $x,y$ and $y,z$, respectively (so the endpoints of $a$ are $x$ and $z$).
We also orient $c$ and $f(c)$, thus assigning signs to the intersection points of $c \cap f(c)$, and so in particular,
to the points $x$, $y$, and $z$.  Two of the signs on $x$, $y$, and $z$ must agree.  If $x$ and $y$, say, have the same
sign, then the curve $a_1 \cup b_1$ is nonseparating since it has geometric intersection number $1$ with the curve $a_1
\cup (c-b_1)$. Therefore, $a_1 \cup b_1$ would be nonperipheral, and we would be done.  Similarly, if $y$ and $z$ have
the same sign, then $a_2 \cup b_2$ is nonseparating and hence nonperipheral.  Therefore, we may assume that the signs
of intersection alternate.

It follows that a regular neighborhood of $a_1 \cup a_2 \cup c = a \cup c$ is as shown in Figure \ref{figure:ac hood},
where we have decomposed $c$ into three arcs $c_1 \cup c_2 \cup c_3$ by the intersection points $x$, $y$, and $z$.
\begin{figure}[ht]
\begin{center}
\psfrag{a}{$a_1$} \psfrag{a'}{$a_2$} \psfrag{c}{$c$} \psfrag{x}{$x$} \psfrag{y}{$y$} \psfrag{z}{$z$}
\psfrag{b_1}{$c_1$} \psfrag{b_2}{$c_2$} \psfrag{b_3}{$c_3$}
\includegraphics[scale=.75]{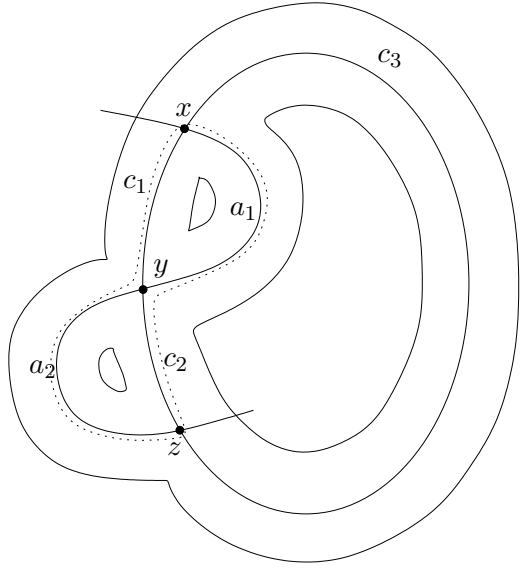}
\caption{Neighborhood of $a_1 \cup a_2 \cup c$ with $c = c_1 \cup c_2 \cup c_3$.  The dotted curve is $d$ (Case 3).} \label{figure:ac hood}
\end{center}
\end{figure}

Each of the arcs $b_1$, $b_2$, and $b$ is made from unions of the three arcs $c_1$, $c_2$, and $c_3$, depending on the
relative lengths of $c_1$, $c_2$, and $c_3$.  There are three cases to consider.\\

\noindent {\bf Case 1.}  $\ell_{q_0}(c_i) \leq \ell_{q_0}(c)/2$ for all $i=1,2,3$.\\

In this case, we have $b_1 = c_1$, $b_2 = c_2$, and $b = c_3$.
Consider the regular neighborhood $N$ of $a \cup c$ shown in Figure
\ref{figure:ac hood}.  We claim that the inclusion of $N$ into $S$
injects on the level of fundamental groups, i.e. $N$ is {\em
incompressible}.  Recall the elementary fact that a subsurface is
incompressible if and only if each of the boundary curves is
homotopically nontrivial (i.e. none of the boundary curves is homotopic
to a point). It follows that $N$ is incompressible since each of the
boundary components is (homotopic to) a union of two segments in $c$ and
$f(c)$ (which were in minimal position), hence homotopically nontrivial.
Note furthermore that $a \cup b$ is not peripheral in $N$, hence cannot be
peripheral in $S$.\\

\noindent {\bf Case 2.} $\ell_{q_0}(c_i) > \ell_{q_0}(c)/2$ for $i=1$ or $i=2$.\\

We consider only the situation $\ell_{q_0}(c_1) > \ell_{q_0}(c)/2$, with the proof for $\ell_{q_0}(c_2) >
\ell_{q_0}(c)/2$ obtained by simply changing the labels.  In this case, we have $b_1 = c_2 \cup c_3$, $b_2 = c_2$, and
$b = c_3$.  We now consider the regular neighborhood $N$ of $a \cup b_1
= a \cup c_2 \cup c_3$, which is a pair of pants.  Note that $a_1 \cup
b_1$, $a_2 \cup b_2$, and $a \cup b$ are all contained in $N$. In fact,
these curves are precisely the three boundary components. Since each of
these curves is a union of two segments in $c$ and $f(c)$, these are
homotopically nontrivial, and so as in Case 1, $N$ is
incompressible. Finally, if all three curves were peripheral, the
complement of $N$ in $S$ would have to be three once-punctured disks,
and hence $S$ would be a thrice-punctured sphere. This is a
contradiction since there are no pseudo-Anosov homeomorphisms of a
thrice-punctured sphere.  Thus, one of the curves must be
nonperipheral.\\

\noindent {\bf Case 3.} $\ell_{q_0}(c_3) > \ell_{q_0}(c)/2$.\\

In this final case, we must have $b_1 = c_1$, $b_2 = c_2$, and $b = c_1 \cup c_2$. Here we let $N$ be the regular
neighborhood of $a \cup b = a_1 \cup a_2 \cup c_1 \cup c_2$.  Again $N$ is a pair of pants, and it contains our three
curves $a_1 \cup b_1$, $a_2 \cup b_2$, and $a \cup b$.  As above $a_1 \cup b_1$ and $a_2 \cup b_2$ are homotopic to two
of the three boundary components, and are both homotopically nontrivial.  If we show that the third boundary component,
$d$ (the dotted curve in Figure~\ref{figure:ac hood}), is homotopically nontrivial, then $a \cup b$, which is an
immersed essential curve in $N$, will be nonperipheral, and this will complete the proof.

If $d$ is homotopically trivial, then it bounds a disk $D$ in $S$.  Since $D$ cannot contain the other two boundary
components of $N$, as these are nontrivial, it follows that $D$ must be ``outside'' of $d$ in Figure \ref{figure:ac
hood}.  We orient $f(c)$ so that it passes through $x$, $y$, and $z$ in that order.  After passing through $z$, $f(c)$
enters $D$.  Since $f(c)$ has no further intersection with $a_1$, $a_2$, $c_1$, and $c_2$ other than the ones shown, it
must cross $c_3$ upon leaving $D$.  But this creates a bigon between $c$ and $f(c)$, contradicting our standing
assumption on minimal position.  It follows that $d$ cannot be homotopically trivial, and hence $N$ is incompressible
and $a \cup b$ is nonperipheral, as required.
\end{proof}

We believe that a much stronger result is true; namely, that dilatations
increase exponentially in the number of
punctures for Brunnian pseudo-Anosov mapping classes.

\begin{conjecture}
There exist constants $A,B>0$ so that
$$L(\brun(S_{g,p})) \geq Ap + B$$
for all $p\geq 1$ and any $g$.
\end{conjecture}


\small

\noindent
Benson Farb:\\
Dept. of Mathematics, University of Chicago\\
5734 University Ave.\\
Chicago, IL 60637\\
E-mail: farb@math.uchicago.edu
\medskip

\noindent
Christopher J. Leininger:\\
Dept. of Mathematics, University of Illinois at Urbana-Champaign \\
273 Altgeld Hall, 1409 W. Green St. \\
Urbana, IL 61802\\
E-mail: clein@math.uiuc.edu
\medskip

\noindent
Dan Margalit: \\
Dept. of Mathematics, University of Utah\\
155 S 1440 East \\
Salt Lake City, UT 84112\\
E-mail: margalit@math.utah.edu


\begin{thebibliography}{ABCDEF}
\small

\bibitem[FLP]{FLP}
{\em Travaux de {T}hurston sur les surfaces}, volume~66 of {\em Ast\'erisque}.
\newblock Soci\'et\'e Math\'ematique de France, Paris, 1979.
\newblock S\'eminaire Orsay, With an English summary.

\bibitem[Ab]{Ab}
William Abikoff.
\newblock {\em The real analytic theory of {T}eichm\"uller space}, volume 820
  of {\em Lecture Notes in Mathematics}.
\newblock Springer, Berlin, 1980.

\bibitem[AY]{AY}
Pierre Arnoux and Jean-Christophe Yoccoz.
\newblock Construction de diff\'eomorphismes pseudo-{A}nosov.
\newblock {\em C. R. Acad. Sci. Paris S\'er. I Math.}, 292(1):75--78, 1981.

\bibitem[BL]{BL}
Hyman Bass and Alexander Lubotzky.
\newblock Linear-central filtrations on groups.
\newblock In {\em The mathematical legacy of Wilhelm Magnus: groups, geometry
  and special functions (Brooklyn, NY, 1992)}, volume 169 of {\em Contemp.
  Math.}, pages 45--98. Amer. Math. Soc., Providence, RI, 1994.

\bibitem[Ba1]{Ba1}
Max Bauer.
\newblock Examples of pseudo-{A}nosov homeomorphisms.
\newblock {\em Trans. Amer. Math. Soc.}, 330(1):333--359, 1992.

\bibitem[Ba2]{Ba2}
Max Bauer.
\newblock An upper bound for the least dilatation.
\newblock {\em Trans. Amer. Math. Soc.}, 330(1):361--370, 1992.

\bibitem[Be]{Be}
Lipman Bers.
\newblock An extremal problem for quasiconformal mappings and a theorem by
  {T}hurston.
\newblock {\em Acta Math.}, 141(1-2):73--98, 1978.

\bibitem[Bi]{Bi}
Joan~S. Birman.
\newblock {\em Braids, links, and mapping class groups}.
\newblock Princeton University Press, Princeton, N.J., 1974.
\newblock Annals of Mathematics Studies, No. 82.

\bibitem[BH]{BH}
Steven~A. Bleiler and Craig~D. Hodgson.
\newblock Spherical space forms and {D}ehn filling.
\newblock {\em Topology}, 35(3):809--833, 1996.

\bibitem[Br]{Br}
Peter Brinkmann.
\newblock A note on pseudo-{A}nosov maps with small growth rate.
\newblock {\em Experiment. Math.}, 13(1):49--53, 2004.

\bibitem[De]{De}
Max Dehn.
\newblock {\em Papers on group theory and topology}.
\newblock Springer-Verlag, New York, 1987.
\newblock Translated from the German and with introductions and an appendix by
  John Stillwell, With an appendix by Otto Schreier.

\bibitem[GT]{GT}
M.~Gromov and W.~Thurston.
\newblock Pinching constants for hyperbolic manifolds.
\newblock {\em Invent. Math.}, 89(1):1--12, 1987.

\bibitem[HS]{HS}
Ji-Young Ham and Won~Taek Song.
\newblock {The minimum dilatation of pseudo-Anosov 5-braids}.
\newblock {Preprint, arXiv:math.GT/0506295}.

\bibitem[H]{H}
W.~J. Harvey.
\newblock Boundary structure of the modular group.
\newblock In {\em Riemann surfaces and related topics: Proceedings of the 1978
  Stony Brook Conference (State Univ. New York, Stony Brook, N.Y., 1978)},
  volume~97 of {\em Ann. of Math. Stud.}, pages 245--251, Princeton, N.J.,
  1981. Princeton Univ. Press.

\bibitem[HK]{HK}
Eriko Hironaka and Eiko Kin.
\newblock A family of pseudo-{A}nosov braids with small dilatation.
\newblock {\em Algebr. Geom. Topol.}, 6:699--738 (electronic), 2006.

\bibitem[Iv1]{Iv1}
Nikolai~V. Ivanov.
\newblock Coefficients of expansion of pseudo-{A}nosov homeomorphisms.
\newblock {\em Zap. Nauchn. Sem. Leningrad. Otdel. Mat. Inst. Steklov. (LOMI)},
  167(Issled. Topol. 6):111--116, 191, 1988.

\bibitem[Iv2]{Iv2}
Nikolai~V. Ivanov.
\newblock {\em Subgroups of {T}eichm\"uller modular groups}, volume 115 of {\em
  Translations of Mathematical Monographs}.
\newblock American Mathematical Society, Providence, RI, 1992.
\newblock Translated from the Russian by E. J. F. Primrose and revised by the
  author.

\bibitem[Jo1]{Jo1}
Dennis Johnson.
\newblock An abelian quotient of the mapping class group {${\cal I}\sb{g}$}.
\newblock {\em Math. Ann.}, 249(3):225--242, 1980.

\bibitem[Jo2]{Jo2}
Dennis Johnson.
\newblock A survey of the {T}orelli group.
\newblock In {\em Low-dimensional topology (San Francisco, Calif., 1981)},
  volume~20 of {\em Contemp. Math.}, pages 165--179. Amer. Math. Soc.,
  Providence, RI, 1983.

\bibitem[Le]{Le}
Christopher~J. Leininger.
\newblock On groups generated by two positive multi-twists: {T}eichm\"uller
  curves and {L}ehmer's number.
\newblock {\em Geom. Topol.}, 8:1301--1359 (electronic), 2004.

\bibitem[MM]{MM}
Howard~A. Masur and Yair~N. Minsky.
\newblock Geometry of the complex of curves. {I}. {H}yperbolicity.
\newblock {\em Invent. Math.}, 138(1):103--149, 1999.

\bibitem[Mc]{Mc}
Curtis~T. McMullen.
\newblock Polynomial invariants for fibered 3-manifolds and {T}eichm\"uller
  geodesics for foliations.
\newblock {\em Ann. Sci. \'Ecole Norm. Sup. (4)}, 33(4):519--560, 2000.

\bibitem[Mk]{Mk}
Hiroyuki Minakawa.
\newblock Examples of pseudo-{A}nosov homeomorphisms with small dilatations.
\newblock {\em J. Math. Sci. Univ. Tokyo}, 13(2):95--111, 2006.

\bibitem[Pe]{Pe}
R.~C. Penner.
\newblock Bounds on least dilatations.
\newblock {\em Proc. Amer. Math. Soc.}, 113(2):443--450, 1991.

\bibitem[So]{So}
Won~Taek Song.
\newblock Upper and lower bounds for the minimal positive entropy of pure
  braids.
\newblock {\em Bull. London Math. Soc.}, 37(2):224--229, 2005.

\bibitem[SKL]{SKL}
Won~Taek Song, Ki~Hyoung Ko, and J{\'e}r{\^o}me~E. Los.
\newblock Entropies of braids.
\newblock {\em J. Knot Theory Ramifications}, 11(4):647--666, 2002.
\newblock Knots 2000 Korea, Vol. 2 (Yongpyong).

\bibitem[Th]{Th}
William~P. Thurston.
\newblock On the geometry and dynamics of diffeomorphisms of surfaces.
\newblock {\em Bull. Amer. Math. Soc. (N.S.)}, 19(2):417--431, 1988.

\bibitem[Wo]{Wo}
Scott Wolpert.
\newblock The length spectra as moduli for compact {R}iemann surfaces.
\newblock {\em Ann. of Math. (2)}, 109(2):323--351, 1979.

\bibitem[Zh]{Zh}
A.~Yu. Zhirov.
\newblock On the minimum dilation of pseudo-{A}nosov diffeomorphisms of a
  double torus.
\newblock {\em Uspekhi Mat. Nauk}, 50(1(301)):197--198, 1995.


\end{thebibliography}
\end{document}